\documentclass{article}
\usepackage{graphicx} 
\usepackage{amsthm, amsmath, amsfonts,mathtools, fullpage, graphicx, wrapfig, tikz, caption, subcaption}
\usepackage{natbib}
\usepackage{amssymb}
\usepackage{verbatim}
\usepackage{dsfont}
\usepackage{relsize}
\usepackage{listings}
\usepackage{booktabs}
\usepackage{xcolor}
\usepackage{hyperref}
\usepackage{graphicx}
\usepackage{multirow}
\usepackage[linesnumbered,ruled,vlined]{algorithm2e}

\newtheorem{theorem}{Theorem}

\newcommand{\E}{\mathbb{E}}
\newcommand{\Prb}{\mathbb{P}}
\newcommand{\cD}{\mathcal{D}}
\newcommand{\cA}{\mathcal{A}}

\title{As Good as it Gets: Bounds for Oracle Time-Varying Treatment Strategies}
\author{Zach Shahn}

\begin{document}
\maketitle

\begin{abstract}
Much causal inference research is focused on methods for optimizing dynamic treatment regimes \citep{murphy2003optimal,robins2004optimal,robins2008estimation,schulte2015q}, which are rules for deciding which treatments should be assigned and when based on evolving history. There is a certain optimism underlying this endeavor that with enough tinkering we might realize consequential improvements. Another strand of research, previously confined to the point exposure setting, considers bounds on how well any individualized treatment rule could possibly do. Here, we extend to the time-varying setting sharp bounds on the performance of an oracle strategy that selects the best treatment regime for each subject based on their unobserved potential outcomes or `response type'. For binary outcomes, the lower bound (assuming higher is better) is simply the expected outcome attained by the optimal treatment regime based on observed history. For continuous outcomes, the lower bound may strictly exceed the maximal observed covariate based value. In the continuous setting, we also consider bounds on the CDF of oracle continuous potential outcomes. 
\end{abstract}

\section{Introduction}

There are countless domains in which treatment decisions are made at multiple time points. In medicine, examples include care for major depressive disorder, epilepsy, HIV, diabetes, sepsis, and many more. Based on evolving patient history, doctors consider: when to initiate treatment; which treatment to try first; when to switch treatment; which treatment to switch to; etc. Given the variety of possible rules (dynamic treatment strategies, or DTRs) for making these decisions, it is hoped that patients might respond favorably to at least one. There exist a range of approaches for estimating the population level optimal DTR \citep{murphy2003optimal,robins2004optimal,schulte2015q} from data satisfying sequential exchangeability and other standard causal assumptions, e.g. from a sequentially randomized trial. The optimal DTR is the rule based on observed history that would lead to the best average outcome if applied to the entire population. However, individuals may still do better under strategies other than the optimal DTR, potentially based on unobserved characteristics. The other side of the coin is that some individuals may not do well under any treatment strategy, including the optimal DTR. The proportions of `doomed' and `treatable' patients \citep{greenland1986identifiability} speak to the potential and limitations of personalization. 

Unfortunately, the `treatable' and `doomed' proportions are not point identified even from a sequentially randomized experiment. In the point exposure setting, some authors have considered point identification under unconventional strong structural assumptions \citep{heckman1997making, wu2025promises}. Others have considered bounds on the potential of personalization \citep{galanter2024can,tian2000probabilities,fan2010sharp,manski1997monotone}. In this note, we develop bounds for the time-varying treatment setting. The bounds are obtained via recursions similar to Bellman equations employed for estimating optimal DTRs. We consider both binary and bounded continuous outcomes. In the continuous setting, we bound both the average utility and the cumulative distribution function (CDF) of outcomes that could be achieved under optimal treatment assignment, extending results from the point exposure setting \citep{fan2010sharp}. 

As a motivating example, consider the ongoing Sequential Multiple Assignment Randomized Trial for Bipolar Depression (SMART-BD) \citep{smartbd2025}.  The trial assigns patients to one of four commonly used treatments and, after six weeks, reassigns `non-remitters' to one of the three treatments they have not yet received. Its primary aims include comparing 12 embedded adaptive interventions and learning whether patient characteristics measured at baseline and six-weeks post-baseline can improve treatment selection. Those questions pertain to how well treatment can be selected using recorded patient information. We instead ask, how much additional remission could conceivably be achieved if the responses to all treatment strategies of each patient were known perfectly?  A small upper bound would suggest that the principal limitation is the treatments themselves rather than the information used to select among them. A large upper bound would leave open a more substantial role for undiscovered effect modifiers.

\section{Setting}

Decisions are made at times $t=0,\ldots,T-1$.  Let $H_t$ denote all information observed before treatment $A_t$ taking values in $\cA$.  The history may include continuous and high-dimensional baseline and time-varying covariates. Let $Y\in\{0,1\}$ be a final outcome, with 1 denoting the preferred outcome. (We will consider bounded continuous $Y$ in section \ref{sec:continuous}.) A dynamic treatment rule $d=(d_0,\ldots,d_{T-1})$ consists of functions
\[
d_t:H_t\longmapsto \cA_t.
\]
Let $\cD$ denote the set of all such rules.

Let $Y(d)$ be the counterfactual outcome under rule $d$.  We consider
\begin{equation}
V^{\mathrm{perfect}}
=\E\left\{\sup_{d\in\cD}Y(d)\right\}.
\label{eq:perfect}
\end{equation}
Because $Y$ is binary, this is the fraction of patients who would succeed under at least one dynamic treatment rule.  It is also the mean outcome that would be obtained if one could select a rule separately for each person with perfect knowledge of all that person's treatment responses.

We make the usual sequential consistency, positivity, and strong exchangeability assumptions. 
\paragraph{Sequential consistency}
For every treatment history $\bar a_t=(a_1,\ldots,a_t)$,
\[
\bar A_t=\bar a_t
\quad\Longrightarrow\quad
H_{t+1}=H_{t+1}(\bar a_t),
\]
and, for a complete treatment history $\bar a_T$,
\[
\bar A_T=\bar a_T
\quad\Longrightarrow\quad
Y=Y(\bar a_T).
\]
Thus, when an individual's observed treatment history agrees with a specified
counterfactual treatment history, the corresponding observed histories and
outcome equal their counterfactual values. This assumption also requires the
treatments to be sufficiently well defined that different versions represented
by the same treatment label can be treated as equivalent.

\paragraph{Sequential positivity}
\[
\Pr(A_t=a\mid H_t=h_t)>0
\]
for every $t$, every $a\in\mathcal{A}_t(h_t)$, and every history $h_t$ that can occur
under a strategy in $\mathcal D$. 

\paragraph{Strong sequential exchangeability}
\[
A_t
\mathrel{\perp\!\!\!\perp}
\mathcal C_t
\mid H_t,
\qquad t=1,\ldots,T,
\]
where
\[
\mathcal C_t
=
\left\{
\begin{array}{l}
H_{t+1}(\bar A_{t-1},a_t),\,
H_{t+2}(\bar A_{t-1},a_t,a_{t+1}),\,\ldots,
Y(\bar A_{t-1},a_t,\ldots,a_T)
\end{array}
\right\}.
\]

Together, these assumptions identify the marginal counterfactual distribution of \(Y(d)\) for every
\(d\in\mathcal D\) through the longitudinal g-formula. They do not,
however, identify the joint distribution of
\(\{Y(d):d\in\mathcal D\}\), and therefore do not point identify
\(V^{\mathrm{perfect}}\). We assume a finest fully randomized causally interpreted structured tree graph model (FFRCISTG) \citep{robins1986new,richardson2013single}.   

\section{Bounds for binary outcomes}
  We begin with backward recursive definitions of quantities that will be used in the bounds. Immediately before the last treatment, let
\[
m_{T-1}(h,a)=\E(Y\mid H_{T-1}=h,A_{T-1}=a).
\]
Set
\begin{align}
q_{T-1}^-(h,a)&=q_{T-1}^+(h,a)=m_{T-1}(h,a),\\
v_{T-1}^-(h)&=\max_{a\in\cA_{T-1}}q_{T-1}^-(h,a),
\label{eq:lastlower}\\
v_{T-1}^+(h)&=\min\left\{1,\sum_{a\in\cA_{T-1}}q_{T-1}^+(h,a)\right\}.
\label{eq:lastupper}
\end{align}
For $t=T-2,\ldots,0$, recursively define
\begin{align}
q_t^-(h,a)
&=\E\left\{v_{t+1}^-(H_{t+1})\mid H_t=h,A_t=a\right\},
\label{eq:qminus}\\
q_t^+(h,a)
&=\E\left\{v_{t+1}^+(H_{t+1})\mid H_t=h,A_t=a\right\},
\label{eq:qplus}\\
v_t^-(h)
&=\max_{a\in\cA_t}q_t^-(h,a),
\label{eq:lower}\\
v_t^+(h)
&=\min\left\{1,\sum_{a\in\cA_t}q_t^+(h,a)\right\}.
\label{eq:upper}
\end{align}
Finally, define
\[
V^-=\E\{v_0^-(H_0)\},
\qquad
V^+=\E\{v_0^+(H_0)\}.
\]
Note that $V^-=\sup_{d\in\cD}\E\{Y(d)\}$ is the value of the optimal DTR based on observed covariates. It should be obvious that the value of the optimal DTR based on observed covariates (i.e. $V^-$) is a lower bound on the value of an oracle DTR, but maybe slightly less so that the bound is sharp as we show in Theorem \ref{theorem1}.

\begin{theorem}[Sharp bounds]\label{theorem1}
Under consistency, sequential exchangeability, positivity, and a FFRCISTG model,
\[
V^{\mathrm{perfect}}\in[V^-,V^+].
\]
The interval is sharp in that both endpoints can be attained by a joint counterfactual distribution consistent with the observed data generating process.
\end{theorem}

The proof of the theorem is in the Appendix, but we sketch it here. Conditional on a history $h$, let $E_a$ be the event that the patient would eventually succeed after choosing treatment $a$ now and then following their optimal oracle regime thereafter. Whatever the association among these events,
\[
\max_a\Prb(E_a\mid h)
\le
\Prb\left(\bigcup_aE_a\mid h\right)
\le
\min\left\{1,\sum_a\Prb(E_a\mid h)\right\}.
\]
The lower bound can be attained by making the success events nested.  The upper bound can be attained by making them disjoint until their total probability reaches one.  The FFRCISTG model permits either.  Applying this argument at each decision time gives \eqref{eq:lower}--\eqref{eq:upper}.

The fraction who would fail under every dynamic rule is $1-V^{\mathrm{perfect}}$. Its sharp bounds are therefore
\begin{equation}
1-V^+\le
\Prb\{Y(d)=0\text{ for every }d\in\cD\}
\le1-V^-.
\label{eq:failure}
\end{equation}
The possible gain from perfect patient-specific response information beyond optimal use of the observed history satisfies the sharp bound
\begin{equation}
0
\le
V^{\mathrm{perfect}}-\sup_{d\in\cD}\E\{Y(d)\}
\le
V^+-V^-.
\label{eq:headroom}
\end{equation}

The bounds can be estimated by methods used to estimate an optimal DTR, e.g. \citep{schulte2015q}.  Following the q-learning approach, estimate the outcome risk $m_{T-1}(h,a)$ at the last decision time. Use the fitted risk to calculate $v_{T-1}^-$ and $v_{T-1}^+$ for each patient. At the preceding time, regress each of these two quantities on the earlier history and treatment to estimate \eqref{eq:qminus} and \eqref{eq:qplus}. Continue backward, and finally average the two fitted baseline values.

\subsection*{A numerical illustration}

The table below contains data from a hypothetical trial with two initial treatments ($A$ and $B$) an early response indicator ($R$), and two second-stage treatments ($C$ and $D$).  The third column gives the probability of early response after the initial treatment, and the last two columns give the conditional probability of final success after each second-stage treatment.

\begin{center}
\begin{tabular}{ccccc}
\hline
Initial treatment & $R$ & $\Prb(R=r\mid\text{initial treatment})$ & $C$ & $D$\\
\hline
$A$ & 1 & .50 & .32 & .04\\
$A$ & 0 & .50 & .02 & .10\\
$B$ & 1 & .40 & .06 & .02\\
$B$ & 0 & .60 & .03 & .01\\
\hline
\end{tabular}\label{tab:numerical1}
\end{center}

After $A$, the lower bounds on the best achievable positive outcome rates are $.32$ for responders and $.10$ for non-responders, while the corresponding upper bounds are $.36$ and $.12$, respectively.  Averaging over $R$ gives
\[
q_0^-(A)=.50(.32)+.50(.10)=.210,
\qquad
q_0^+(A)=.50(.36)+.50(.12)=.240.
\]
The same calculations after $B$ give $q_0^-(B)=.042$ and $q_0^+(B)=.056$.  Therefore,
\[
V^-=\max(.210,.042)=.210,
\qquad
V^+=.240+.056=.296.
\]
The fraction who would succeed under at least one treatment rule is sharply bounded by $[.210,.296]$.  The optimal implementable rule starts with $A$, then assigns $C$ to early responders and $D$ to non-responders and has success probability $.210$. Perfect knowledge of individual treatment responses could improve the success rate by at most $.086$. The interval is well below one and fairly narrow because all alternatives to the optimal rule have small probability of success.

\section{Bounds for bounded continuous outcomes}\label{sec:continuous}

Suppose that the outcome under each strategy is bounded, with
$Y(d)\in[0,1]$, and define the best outcome available to an oracle as
\[
M=\sup_{d\in\mathcal D}Y(d).
\]
The distribution of $M$ can be bounded by applying the binary-outcome
recursion separately at each outcome threshold. For $x\in[0,1]$, define
\[
Y_x(d)=\mathbb{I}\{Y(d)>x\}.
\]
Because
\[
\{M>x\}
=
\bigcup_{d\in\mathcal D}\{Y(d)>x\},
\]
the probability $\Pr(M>x)$ is precisely the oracle value for the binary
outcome $Y_x(d)$.

At threshold $x$, replace the terminal success probabilities in the
recursion by the corresponding conditional `survival' probabilities,
\[
m_x(h,a)
=
\Pr(Y>x\mid H_T=h,A_T=a),
\]
and otherwise apply the same backward recursion. Let the resulting
endpoints be $V_x^-$ and $V_x^+$. Then
\[
V_x^-
\leq
\Pr(M>x)
\leq
V_x^+,
\qquad x\in[0,1].
\]
Equivalently, if $F_M(x)=\Pr(M\leq x)$, then
\[
1-V_x^+
\leq
F_M(x)
\leq
1-V_x^-.
\]
Thus the method bounds the entire distribution of the best outcome that
could be achieved if the best strategy for each individual were known.
For a clinically meaningful threshold $x$, the bounds describe the
possible proportion of individuals for whom at least one dynamic
strategy would produce an outcome exceeding that threshold. They can
also be inverted to bound quantiles of the oracle outcome distribution.

It turns out these are sharp bounds on the complete distribution, not just (as it might seem) merely
pointwise sharp bounds considered separately at each threshold. To see
why, first define the counterfactual outcome $M_a(h_t)=sup_{d_{t+1:T}}Y(h_t,a,d_{t+1:T})$ at a given history $h_t$ through time $t$ under treatment $a$ at $t$ followed by the oracle optimal regime thereafter. For notational simplicity, we will write this as $M_a$ and suppress the conditioned-upon history. Write the marginal
distribution function of $M_a$ as $F_a$. To attain the lower endpoint, take a single $U\sim Uniform(0,1)$ and set $M_a=F_a^{-1}(U)$ for every action $a$. This pairs the same outcome quantiles across actions. Individuals with relatively favorable outcomes under one action also have relatively favorable outcomes under the others. Under this distribution,
\[
\Pr\left(\max_a M_a\leq x\right)
=
\min_a F_a(x),
\]
or, equivalently,
\[
\Pr\left(\max_a M_a>x\right)
=
\max_a\{1-F_a(x)\}.
\]
For the upper bound, \citet{lai1976maximally} show that for any finite collection of continuous marginal distribution functions $F_a$, there exists a single joint distribution with these marginals under which
\[
\Pr\left(\max_a M_a>x\right)
=
\min\left\{1,\sum_a\{1-F_a(x)\}\right\}
\]
simultaneously for every $x$. Thus, the pointwise union bounds are compatible with a single joint distribution of the continuous outcomes and together form an attainable bound on the entire CDF. This construction can trivially be propagated backward through the treatment tree.

The CDF bounds also yield sharp bounds on the oracle mean. Since
$M\in[0,1]$,
\[
V^{\mathrm{perfect}}
=
\mathbb{E}(M)
=
\int_0^1 \Pr(M>x)\,dx,
\]
and hence
\[
\int_0^1 V_x^-\,dx
\leq
V^{\mathrm{perfect}}
\leq
\int_0^1 V_x^+\,dx.
\]
Unlike the binary case, the lower endpoint is not necessarily the value of any implementable
strategy. At each threshold, $V_x^-$ is the largest probability of
exceeding that threshold attainable by an implementable strategy, but
the strategy attaining this value may differ across thresholds.
Consequently,
\[
\int_0^1 V_x^-\,dx
\geq
\sup_{d\in\mathcal D}\mathbb{E}\{Y(d)\},
\]
while still providing a valid lower bound on the oracle mean. 

In practice, bounds can be estimated
over a grid of $x$ and the mean bounds obtained by numerical
integration. Sampling error may cause the estimated bounds to violate
monotonicity, in which case they can be projected onto
monotone functions. 

\subsection*{A numerical illustration}

Consider another stylized two-stage trial. At baseline, participants are randomized
to $A_1\in\{0,1\}$. Under either initial treatment, the probability of an
early response $R=1$ is $0.5$. Responders continue their initial treatment,
and nonresponders are randomized between two second-stage treatments
$A_2\in\{0,1\}$.

The final outcome lies in $[0,1]$. Within each treatment branch, suppose
\[
Y\mid A_1=a_1,R=r,A_2=a_2
\sim
\operatorname{Beta}\{50\mu,50(1-\mu)\},
\]
where $\mu$ is the corresponding entry below:
\[
\begin{array}{c c c c}
\toprule
A_1 & R & \mu(A_2=0) & \mu(A_2=1)\\
\midrule
0 & 0 & 0.46 & 0.50\\
0 & 1 & \multicolumn{2}{c}{0.52\text{ (continue)}}\\
1 & 0 & 0.35 & 0.31\\
1 & 1 & \multicolumn{2}{c}{0.65\text{ (continue)}}\\
\bottomrule
\end{array}
\]
The beta distributions are continuous, have means $\mu$, and are moderately
concentrated around those means.

Let $S_{a_1,r,a_2}(x)$ denote the conditional probability of exceeding
threshold $x$. Among nonresponders, the second-stage recursion gives
\[
v_{a_1,0}^-(x)
=
\max_{a_2}S_{a_1,0,a_2}(x)
\]
and
\[
v_{a_1,0}^+(x)
=
\min\left\{
1,\sum_{a_2}S_{a_1,0,a_2}(x)
\right\}.
\]
For responders there is no second-stage choice, so
\[
v_{a_1,1}^-(x)
=
v_{a_1,1}^+(x)
=
S_{a_1,1}(x).
\]
Averaging over response gives
\[
q_{a_1}^{\pm}(x)
=
0.5v_{a_1,0}^{\pm}(x)
+
0.5v_{a_1,1}^{\pm}(x),
\]
and the baseline recursion gives
\[
V_x^-=\max_{a_1}q_{a_1}^-(x),
\qquad
V_x^+=\min\{1,q_0^+(x)+q_1^+(x)\}.
\]

The resulting exceedance-probability and CDF bounds at selected thresholds
are
\[
\begin{array}{c c c c c}
\toprule
x & V_x^- & V_x^+ & 1-V_x^+ & 1-V_x^-\\
\midrule
0.40 & 0.939 & 1.000 & 0.000 & 0.061\\
0.50 & 0.556 & 1.000 & 0.000 & 0.444\\
0.55 & 0.465 & 0.804 & 0.196 & 0.535\\
0.60 & 0.387 & 0.502 & 0.498 & 0.613\\
0.70 & 0.117 & 0.120 & 0.880 & 0.883\\
\bottomrule
\end{array}
\]
where $[1-V_x^+,1-V_x^-]$ bounds $\Pr(M\leq x)$.

Numerically integrating the exceedance bounds gives
\[
0.552
=
\int_0^1V_x^-\,dx
\leq
V^{\mathrm{perfect}}
\leq
\int_0^1V_x^+\,dx
=
0.613.
\]
Thus the unidentified association among counterfactual treatment branches
creates an oracle-mean interval only $0.061$ units wide, and its upper
endpoint remains far below one.

The best implementable mean is smaller. Starting with $A_1=0$ and assigning
nonresponders to $A_2=1$ gives
\[
0.5(0.52)+0.5(0.50)=0.510,
\]
whereas starting with $A_1=1$ and assigning nonresponders to $A_2=0$ gives
\[
0.5(0.65)+0.5(0.35)=0.500.
\]
The optimal implementable value is therefore $0.510$. By contrast, the
integrated lower oracle bound is $0.552$. This occurs because $A_1=0$
maximizes the probability of exceeding lower thresholds, whereas $A_1=1$
maximizes it above approximately $x=0.511$. The integrated lower bound takes
the upper envelope of these two exceedance curves, even though no single
implementable strategy attains that envelope at every threshold.

Relative to the best implementable strategy, the possible value of perfect
individual-level treatment selection is consequently bounded by
\[
0.552-0.510
\leq
V^{\mathrm{perfect}}-V^{\mathrm{opt}}
\leq
0.613-0.510,
\]
or
\[
0.042
\leq
V^{\mathrm{perfect}}-V^{\mathrm{opt}}
\leq
0.103.
\]

\section{Discussion}

The bounds presented quantify the greatest improvement that perfect knowledge of
individual treatment responses could provide beyond strategies based only
on observed history. They can therefore be relevant to whether effort should be devoted to discovery of new biomarkers to guide treatment decisions, or whether entirely new treatments are needed. 

Numerical examples were designed such that the bounds were informative. However, they may become less so as the treatment horizon or the
number of available viable actions increases. Repeated applications of the
union bound can cause the upper endpoint to reach one, particularly when
there are many moderately successful strategies and no strategy
clearly dominates the others. Longer treatment sequences also create
statistical difficulties because later conditional distributions must be
estimated from increasingly sparse histories. The bounds may
nevertheless remain informative when success probabilities are small,
when the treatment menu is limited, or when one treatment option is
consistently superior.

Statistical inference for the bound endpoints is nonstandard because the
recursion repeatedly applies maxima, minima, and truncation. The ordinary
bootstrap need not be valid. Methods around these difficulties have been considered in the point exposure setting \citep{huang2017inequality}, and extending them to time-varying treatments is a direction for future work.

\bibliography{bibliography}

\section*{Appendix}
\subsection*{Appendix A: Proof of the theorem}

We prove sharpness by backward induction.  The argument is conditional on each realized history. At the last decision time, the success probability under each treatment $a$ is $m_{T-1}(h,a)$. The probability of success under at least one treatment is the probability of a union of events having these marginals. Its sharp lower and upper bounds are \eqref{eq:lastlower} and \eqref{eq:lastupper}. For the lower endpoint, place all success events on a common uniform random variable so that they are nested.  For the upper endpoint, place intervals of the required lengths consecutively in the unit interval so that their union has probability equal to the sum of their lengths capped at 1.

Now suppose that $v_{t+1}^-(h'_{t+1})$ and $v_{t+1}^+(h'_{t+1})$ are sharp conditional bounds after every possible history $h'_{t+1}$ through time $t+1$.  Conditional constructions attaining those endpoints can be selected separately within each value and then averaged over the identified distribution of $H_{t+1}(a)$, the counterfactual value of $H_{t+1}$ under treatment $a$ at time $t$ and observed treatment trajectory consistent with $H_t$. This gives the attainable probabilities $q_t^-(h,a)$ and $q_t^+(h,a)$ in \eqref{eq:qminus}--\eqref{eq:qplus}.

Fix a decision time $t$ and history $h$. For each available action $a$, let $\mathcal{B}_a$ denote the entire counterfactual subtree that follows if action $a$ is taken at $h$, including all subsequent counterfactual variables in that branch. Let $E_a$ be the event that at least one path within $\mathcal{B}_a$ ultimately succeeds. By the induction hypothesis, we can choose a compatible distribution within each branch such that $\Pr(E_a \mid H_t=h)=q_t^-(h,a)
$
or $\Pr(E_a \mid H_t=h)=q_t^+(h,a)$. These branch-specific distributions are consistent with all identified quantities.

What remains unspecified is how the different branch objects
$\{\mathcal{B}_a : a\in\mathcal{A}\}$ are associated with one another. Because they correspond to mutually incompatible actions at the same history, there are no cross-world restrictions on their association under a FFRCISTG. We may therefore combine them so that their success events have any joint distribution compatible with their marginal probabilities.

For the lower endpoint, take a common $U\sim\operatorname{Uniform}(0,1)$ and set
\[
\mathbb{I}(E_a)=\mathbb{I}\{U\leq q_t^-(h,a)\}.
\]
The events are then nested, so
\[
\Pr\left(\bigcup_a E_a \,\middle|\, H_t=h\right)
=
\max_a q_t^-(h,a)
=
v_t^-(h).
\]

For the upper endpoint, choose sets $J_a\subseteq[0,1]$ having lengths $q_t^+(h,a)$, arranged consecutively, and set
\[
\mathbb{I}(E_a)=\mathbb{I}\{U\in J_a\}.
\]
Their union has the sum of their lengths capped at one, and hence
\[
\Pr\left(\bigcup_a E_a \,\middle|\, H_t=h\right)
=
\min\left\{1,\sum_a q_t^+(h,a)\right\}
=
v_t^+(h).
\]

These constructions can be made without changing the distribution within any branch. Conditional on the assigned value of $\mathbb{I}(E_a)$, draw $\mathcal{B}_a$ from its previously constructed distribution conditional on $E_a$ or $E_a^c$, respectively. Because $\mathbb{I}(E_a)$ has the correct marginal probability, mixing these two conditional distributions recovers the original marginal distribution of $\mathcal{B}_a$. Thus only the association between mutually incompatible branches has been changed, which is permissible. Applying this construction successively from the final decision backward to baseline yields compatible joint counterfactual distributions attaining $V^-$ and $V^+$. This proves sharpness of the interval.

\end{document}